\documentclass[12pt,oneside,a4paper,leqno]{article}

\usepackage[all]{xy}
\CompileMatrices

\newdir{ >}{!/10pt/@{ }*@{>}}
\newdir{< }{!/10pt/@{ }*@{<}}
\usepackage{geometry}
\usepackage{graphicx}
\usepackage{amsfonts,amssymb}
\usepackage[centertags]{amsmath}

\usepackage{amsthm}
\newcounter{theorems}

\theoremstyle{plain}
\newtheorem{theoAA}[theorems]{Theorem}

\swapnumbers
\newcounter{lemma}
\numberwithin{equation}{section}

\newtheoremstyle{par}%
     {\topsep}%
     {\topsep}%
     {\itshape}%
     {}%
     {\bfseries}%
     {}%
     {.5em}%
     {}%

\newtheoremstyle{parrm}%
     {\topsep}%
     {\topsep}%
     {\normalfont}%
     {}%
     {\itshape}%
     {}%
     {.5em}%
     {}%
\theoremstyle{plain}
\numberwithin{equation}{section}

\theoremstyle{definition}
\newtheorem{defi}[equation]{Definition}
\newtheorem{example}[equation]{Example}

\theoremstyle{remark}
\newtheorem{remark}[equation]{Remark}

\theoremstyle{par}
\newtheorem{nr}[equation]{}
\newtheorem{lemma}[equation]{}
\newtheorem{propo}[equation]{}

\theoremstyle{parrm}

\makeatletter
\def\tagform@#1{\maketag@@@{\ignorespaces#1\unskip\@@italiccorr}}

\makeatother
\usepackage{paralist}
\usepackage{latexsym,amssymb}
\newcommand{\minus}{\smallsetminus}

\newcommand{\ze}{\mathbb{Z}}
\newcommand{\st}{\ \mathrm{|} \ }
\newcommand{\from}{\colon}

\newcommand{\action}{\mathcal{A}}

\newcommand{\kk}{\mathbf{k}}

\newcommand{\nn}{\mathbf{n}}

\newcommand{\dt}{\ dt}
\newcommand{\menouno}{{-1}}
\newcommand{\ZZ}{\mathbb{Z}}

\newcommand{\bfomega}{{\bf\omega}}
\newcommand{\TT}{\mathbb{T}}
\newcommand{\RR}{\mathbb{R}}
\newcommand{\XX}{{X}}
\newcommand{\Iso}{\mathop{\mathrm{Iso}}\nolimits}
\newcommand{\IsoR}{{\Iso(\RR)}}
\newcommand{\IsostarR}{{\Iso_*(\RR)}}
\renewcommand{\SS}{\mathbb{S}}

\newif\ifpdf%
\ifx\pdfoutput\undefined%
\pdffalse%
\else%
\pdfoutput=1%
\pdftrue%
\fi%
\newcommand{\tG}{{\widetilde G}}
\newcommand{\Krh}{Krh}
\newcommand{\hKrh}{\widehat{Krh}}
\begin{document}
\pagenumbering{arabic}
\title{%
Transitive decomposition of symmetry groups for the $n$-body problem
}

\author{Davide L.~Ferrario}
\date{%
\today}
\maketitle

\begin{abstract}
Periodic and quasi-periodic orbits of the $n$-body problem 
are critical points of the action functional constrained to the Sobolev
space of symmetric loops. Variational methods yield collisionless
orbits provided the group of symmetries fulfills certain conditions 
(such as the \emph{rotating circle property}). Here we generalize 
such conditions to more general group types  and  show
how to constructively classify all groups satisfying such hypothesis,
by a decomposition into irreducible transitive components. As examples we 
show approximate trajectories of 
some of the  resulting symmetric minimizers.
\end{abstract}
\section{Introduction and main results}
\label{sec:intro}
Periodic and quasi-periodic for the $n$-body problem have received
much of 
attention over the last years, also because of the success of 
variational and topological methods. 
The starting point can be traced back to  the 
nonlinear analysis  works
of 
A. Ambrosetti, A. Bahri, V. Coti-Zelati, P. Majer, J. Mawin, P.H. Rabinowitz,
E. Serra 
and S. Terracini (among others) issued around 1990 
\cite{%
bara91,%
mawi89,%
ambrosetti92,%
ambrocoti,%
ambrosetticoti,%
mate95,%
mate93,%
serter2%
};  methods were developed that could deal with singular potential and particular
symmetry groups of the functional. 
For other approaches one can see also
 I. Stewart \cite{stewart} and  C. Moore \cite{MR1220207}.
The next new wave of results has followed the 
remarkable A. Chenciner and R. Montgomery's  proof of  the ``figure-eight'' 
periodic solution of the three-body problem in the case of equal masses
\cite{monchen}, where collisions and singularities were excluded by
the computation of the action level on test curves and a non-commutative 
finite group of symmetries was taken as the constraint for a global
equivariant variational approach. In order to generalize the equivariant 
variational method  (that is, to restrict the action functional to 
the space of equivariant loops)
to a new range 
of applicability S. Terracini and the author in \cite{FT2003}
could make use of C. Marchal's averaging idea \cite{marchal} and 
prove that local minimizers of the action functional are collisionless,
provided an algebraic condition on the symmetry group  
(termed the \emph{rotating circle property}) holds.
Meanwhile, symmetry groups and various approaches to level estimates 
or local variations have been found, together with the corresponding
symmetric minimizers, and published by many authors (see 
for example
\cite{BT2004,zz,dd},
\cite{MR2020261,MR2032484},
\cite{chenven},
\cite{MR1916506},
\cite{MR2146342},
\cite{MR1963764},
\cite{MR1919782,MR2134901},
and \cite{MR2012847}.
The aim of this article is to provide a unified framework for the 
construction and classification 
problem and at the same time to extend the application range 
of the averaging and blow-up techniques.
More precisely, when dealing with the problem of classifying in a constructive way 
finite symmetry groups
for the $n$-body problem one has to face three  issues.
First, it is of course preferrable to have an equivalence relation  defined
between groups, which  rules out differences 
thought as non--substantial. Second, one has to find a  suitable
decomposition of a symmetry group into a  sum of 
(something like) \emph{irreducible} components. The way of decomposing
things depends upon the context. In our settings we could choose 
an orthogonal  representations decomposition (as direct sum of $G$-modules)
or as permutations decomposition, or a mixture of the two. Third, it would be
interesting to deduce from the irreducible components, from their
(algebraic and combinatorial) properties, some consequente 
properties
of action-minimizing periodic orbits (like being collisionless,
existence, being homographic or non-homographic, and the like).
The purpose of the paper is to give a procedure for 
constructing all symmetry groups 
of the $n$-body problem (in three-dimensional
space, but of course the planar case can be done as a particular
case) according to these three options, with the main focus
on the existence of periodic or quasi-periodic non-colliding solutions. 
The main result can be used to 
list groups that might be considered as the
elementary {building blocks} for generic symmetry groups
yielding collisionless minimizers. 

The first reduction  will be obtained by defining  the \emph{cover}
of a symmetry group (that is, the group acting on the time line
instead of the time circle) and considering equivalent groups 
with the same equivariant periodic
trajectories (up to repeating loops). Also, it is possible 
to consider equivalent those symmetry groups that differ by a change 
in the action functional (the angular speed).
Using this simple \emph{escamotage} it is possible to dramatically reduce 
the cardinality of the symmetry groups and to 
deal with a finite number of (numerable) families  of groups for every
$n$. The next step is to exploit the fact that any finite permutation representation
can be decomposed \emph{\`a la Burnside} into  the disjoint
sum of \emph{transitive} (or, equivalently,
\emph{homogeneous}) 
permutation representations. This decomposition requires the definition
of a suitably crafted sum of Lagrangean symmetry groups,
which will be written in term of $\Krh$ and $\hKrh$ data  (to be defined
later) yielded by
the group. The transitive decomposition allows to state
the main result, which  can be written as follows. 
Definition and notation of course
refer to the body of the paper. 

\begin{theoAA}
\label{theo:main}
Let $G$ be a symmetry group with a
colliding $G$-symmetric Lagrangean local minimizer.
If $G_*$  is the $\TT$-isotropy group of  
the colliding time restricted to the index subset of colliding bodies, 
then $G_*$ cannot act trivially 
on the index set; if the permutation  isotropy of  a transitive component of  $G_*$ 
is trivial, than the image of $G_*$ in $O(3)$ cannot be 
one of 
the following:
$I$, 
$C_p$ (for $p\geq 1$), 
$D_p$ (for $p\geq 2$), 
$T$, 
$O$, 
$Y$,
$P'_{2p}$,
$C_{ph}$.
\end{theoAA}

This result allows to clarify and to extend the above-mentioned 
\emph{rotating circle
property};
in the proof we show how with a simple application of the averaging 
Marchal technique on space equivariant spheres  one can deduce 
that for the group actions listed in the statement minimizers are
collisionless. It is also worth mentioning that the transitive decomposition
approach has two interesting consequences: from one hand it is possible
to determine whether the hypothesis of Theorem \ref{theo:main} is fulfilled
simply by computing the space-representations of the transitive 
decomposition of the maximal $\TT$-isotropy subgroups of the 
symmetry group (thus making the task of finding
rotating circles unnecessary); on the other hand 
a machinery for finding examples of symmetry groups can be 
significantly improved by allowing the construction of 
actions using smaller and combinatorial components.
Even if feasible, a complete classification 
of all symmetry groups
satistying the hypotheses which imply collision-less minimizers  and
coercivity
would just result into an unreasonably long unreadable list.
We decided to formulate only the method that can be used for such 
generation, leaving a few examples   in the last section 
to illustrate it in simple cases.
Therefore the paper is basically organized as  a multi-step proof of Theorem 
\ref{theo:main}, together with the introduction and explanation 
of the necessary preliminaries, results and definitions. 
In sections 2 and 3 we will review the
notation about Euclidean symmetry groups
and the main properties of Lagrangean symmetry groups. 
In section 4 the definition of transitive decomposition and disjoint
sum of symmetry groups is carried out: this is the one of the main
step in the construction process. Furthermore, in section 5 a simple
proof allows to extend the averaging technique to all orientation-preserving
finite isotropy groups. Together with the rotating circle
property and the classification of finite
subgroups of $SO(3)$ this will yield the method on avoiding collisions.
The analysis of possible transitive component is then carried out in 
section 6, according to the previous definition and results.
At the end, in section 7 the few examples mentioned above  are shown,
together with pictures of the corresponding approximate minimizers.

\section{Preliminaries and notation}
\label{sec:prelim}
We denote by $O(d)$  
the orthogonal group in dimension $d$,
that is, the group of $d\times d$ orthogonal matrices 
over the real field. The symbol $\Sigma_n$ denotes the permutation
group on $n$ elements
$\{1, \dots, n\}$.
Space isometries are named \emph{rotation},
\emph{reflection}, \emph{central inversion}
and \emph{rotatory reflection} (actually a central
inversion is a particular rotatory reflection).
We recall, following the terminology and notation of 
\cite{coxeter} (page 99, pages 270--277)
and
\cite{MR1079192} (appendix A, pages 351--367; see also
\cite{bartsch} and \cite{miller}),
that the non-trivial finite subgroups of $SO(3)$  are
the following: 
$C_p$ (the cyclic group generated by a rotation of order $p$, 
for $p\geq 2$, with a single $p$-gonal axis), 
$D_p$ (the dihedral group of order $2p$, with $p$ horizontal digonal axes 
and a vertical $p$-gonal axis, with $p\geq 2$),
$T\cong A_4$ (the tetrahedral group  of order $12$,
with $4$ trigonal axes and $3$ mutually orthogonal digonal axes),
$O\cong S_4$ (the octahedral group of order $24$,  with $4$ trigonal axes, 
the same as $T$,  and $3$ mutually
orthogonal tetragonal axes;
it is isomorphic to the orientation-preserving 
symmetry group of the cube and contains the tetrahedral group as a normal
subgroup of index $2$) and 
$Y\cong A_5$ (the icosahedral group of order $60$,
with $6$ pentagonal axes, $10$ trigonal axes and $15$ digonal axes).
The dihedral group $D_2$ is a normal subgroup of $T$ of index $3$.

The finite subgroups of $O(3)$ are index $2$ extensions  
of the groups listed
above. Let $I$ denote the group generated
by the central inversion $\menouno\in O(3)$. Since 
$O(3) =  I \times SO(3)$ and $I$ is the center of $O(3)$,
finite groups containing the central inversion are 
$I$, $I\times C_p$, $I\times D_p$, 
$I\times T$, $I\times O$ and $I\times Y$.

The remaining \emph{mixed groups}
are those not containing the central inversion: 
$C_{2p}C_p$ (of order $2p$),
$D_p C_p$ (of order $2p$, it is a Coxeter group, i.e. generated by plane 
reflections; it is
the full symmetry group of a $p$-gonal pyramid),
$D_{2p} D_p$ (of order $4p$; it is a Coxeter group if $p$ is odd, full symmetry
group of a $p$-gonal prism or a $p$-gonal dipyramid)
and $S_4A_4 = OT$ (of order $24$, it is a Coxeter group: the full symmetry
group of a tetrahedron). One word about notation:
mixed groups are denoted by a pair $GH$, where $G$ is a finite rotation
group  of table \ref{tab:finitegroups1}, 
which turns out to be isomorphic to the group under observation
but not conjugated to it, and $H$ a subgroup of index $2$ in $G$. Given such a pair,
a group not containing $I$ 
is obtained as the union (of sets) $H \cup (\menouno(G\minus H))$.
Let $\zeta_p$ and $\kappa$ be the rotations 
\[
\zeta_p   = 
\begin{bmatrix}
\cos 2\pi /p & -\sin 2\pi/p & 0 \\
\sin 2\pi /p & \cos 2\pi/p & 0 \\
0  & 0  & 1 \\
\end{bmatrix}, 
\ \ 
\kappa = 
\begin{bmatrix}
1 & 0  & 0 \\
0  & -1  & 0 \\
0  & 0  & -1 \\
\end{bmatrix}, 
\]
and, if $\varphi=(\sqrt{5}+1)/2$ 
denotes the golden ratio, let  $\pi_3$ and 
$\pi_3'$ be the rotations defined by the following matrices.
\[
\pi_3   =
\begin{bmatrix}
0  & 1  & 0 \\
0  & 0  & 1 \\
1 & 0  & 0 \\
\end{bmatrix}, 
\ \
\pi'_3 = 
\begin{bmatrix}
\varphi/2 & (1-\varphi)/2 & 1/2 \\
(\varphi - 1)/2 & -1/2 & - \varphi/2 \\
1/2 & \varphi/2 & ( 1-\varphi )/2 \\
\end{bmatrix}
\]
Then the generators and normalizers of 
finite subgroups of $SO(3)$ are listed in table \ref{tab:finitegroups1}.
For more data on the icosahedral group, see also \cite{MR0080930}.
\begin{table}
\caption[Finite subgroups of $SO(3)$]{Finite subgroups of $SO(3)$, their normalizers in $SO(3)$ and generators (the generators of  the normalizer
are obtaining adding the generator of the fourth column to the generators of the
second column)}
\label{tab:finitegroups1}
\begin{tabular}{ll|llll}
\emph{Name} & \emph{Symbol} & \emph{Order} &  \emph{Gen.} & \emph{$N_{SO(3)}G$} & \emph{Gen.} \\
\hline
Rotation Cyclic   &  $C_p$ &  $p\geq 2$ & $\zeta_p$ &  $O(2)$ 
& $\zeta_*$, $\kappa$ \\
Rotation Three Axes    &  $D_2$ & $4$ &  $\zeta_2$, $\kappa$  & $O$    & $\zeta_4$, $\pi_3$ \\ 
Rotation Dihedral &  $D_p$ & $2p\geq 6$ &  $\zeta_p$, $\kappa$   &  $D_{2p}$ & $\zeta_{2p}$ \\
Rotation Tetrahedral &  $T\cong A_4$ & $12$ &   $\zeta_2$, $\pi_3$ & $O$ & $\zeta_4$ \\ 
Rotation Octahedral & $O\cong S_4$ & $24$  & $\zeta_4$, $\pi_3$ & $O$ &  \\
Rotation Icosahedral &  $Y\cong A_5$  & $60$  & $\pi_3$, $\pi_3'$  & $Y$ & \\
\end{tabular}
\end{table}

\begin{table}
\caption[Finite subgroups of $O(3)$]{Finite subgroups of $O(3)$ containing
the central inversion} 
\label{tab:finitegroups2}
\begin{tabular}{ll|llll}
\emph{Name} & \emph{Symbol} & \emph{Order} &  \emph{Gen.} & \emph{$N_{O(3)}G$} & \emph{Gen.} \\
\hline
Central Inversion &  $I$ &  $2$ & $\menouno$ &  $O(3)$  & \\
Prism/Antiprism &  $I\times C_p$ &  $2p\geq 4$ & $\menouno$, $\zeta_p$ &  $I \times O(2)$ 
& $\zeta_*$, $\kappa$ \\
Three planes &  $I\times D_2$ & $8$ &  $\menouno$, $\zeta_2$, $\kappa$  & $I\times O$    & $\zeta_4$, $\pi_3$ \\ 
&  $I \times D_p$ & $2p\geq 6$ &  $\menouno$, $\zeta_p$, $\kappa$   &  $I\times D_{2p}$ & $\zeta_{2p}$ \\
&  $I \times T$ & $24$ &   $\menouno$, $\zeta_2$, $\pi_3$ & $I\times O$ & $\zeta_4$ \\ 
Full octahedron  & $I\times O$ & $48$  & $\menouno$, $\zeta_4$, $\pi_3$ & $I\times O$ &  \\
Full icosahedron &  $I \times Y$  & $120$  & $\menouno$, $\pi_3$, $\pi_3'$  & $I\times Y$ & \\
\end{tabular}
\end{table}

\begin{table}
\caption[Finite subgroups of $O(3)$]{Finite subgroups of $O(3)$ 
of mixed type,
their normalizers
and generators}
\label{tab:finitegroups3}
\begin{tabular}{ll|llll}
\emph{Name} & \emph{Symbol} & \emph{Order} &  \emph{Gen.} & \emph{$N_{O(3)}G$} & \emph{Gen.} \\
\hline
Prism/Antisprism & $C_{2p}C_p$ & $2p\geq 2$ & $-\zeta_{2p}$  &  $I\times O(2)$ &  $\zeta_*$,$\menouno$ \\
Reflections Dihedral & $D_pC_p$ & $2p\geq 4$ & $\zeta_p$, $-\kappa$ &  $I\times D_{2p}$ & $\menouno$ \\
 & $D_{2p}D_p$ & $4p\geq 4$ & $\zeta_p$, $\kappa$, $-\zeta_{2p}$ & $I\times D_{2p}$ & $\menouno$ \\
Full tetrahedron & $OT$ & $24$ &  $\zeta_2$,$\pi_3$,$-\zeta_4$ & $I\times O$ &  $\menouno$\\
\end{tabular}
\end{table}

Note that other symbols might be used:
$O^-=\bar{O} = OT$,
$D^d_{2p} = D_{2p}D_p$,
$D^z_p = D_pC_p$,
$Z_{2p}^- = \bar{Z}_{2p} = C_{2p}C_p$,
$\ZZ_2^c = I$,
$Z_p = C_p$,
$I=Y$ (here there is a notation clash with $I=\langle \menouno\rangle$);
the Sch\"onflies notation for crystallographic
point groups (or the equivalent Hermann-Mauguin notation)
is also another option:
for example, $T_d=OT$,
$T_h=I\times T$,
$O_h=I\times O$,
$Y_h=I\times Y$ or
$D_pC_p=C_{pv}$.
Groups generated by reflections (that is, Coxeter groups) are
$D_pC_p$  (with $p\geq 1$),  
$D_{2p}D_p$ (with $p$ odd), 
$I \times D_{p}$ (with $p$ even), 
$OT$, 
$I \times O$, 
$I \times Y$.

Finally, note that the $G$-orbit of a point in general position in 
$\RR^3$ is a regular $p$-agon  for $G=C_p$ 
but it is not a regular polygon for $G=D_p$ or
if $G$ is a polyhedral group (full or rotation).
For the groups $I\times C_p$ and 
$C_{2p}C_p$,
 the $G$-orbit of a point (in general position
in $\RR^3$) is the set of vertices of a \emph{prism}
in $p$ is even and $G=I\times C_p$ 
or if $p$ is odd and $G=C_{2p}C_p$.
It is an \emph{antiprism} (also known as \emph{twisted prism}) 
if $p$ is odd
and $G=I\times C_p$ or $p$ is even and $G=C_{2p}C_p$.
Therefore such groups might be called 
prism/antiprism groups correspondingly. 
In the Schoenflies notation
the antiprism group of order $2p$ is denoted by $S_{2p}$
and the prism group of order $2p$ by $C_{ph}$.  
To avoid possible confusion, we define for $p\geq 1$ the antiprism
group $S_{2p}$ also as
\[
P'_{2p} = 
\begin{cases}
	I\times C_p & \text{if $p$ is odd} \\
        C_{2p}C_p & \text{if $p$ is even.} 
\end{cases}
\]   
It is a cyclic group generated by a \emph{rotatory reflection} of order $2p$.
The prism group on the other hand is defined  for $p\geq 1$  as 
\[
C_{ph} = 
\begin{cases}
 I\times C_p & \text{if $p$ is even} \\
        C_{2p}C_p & \text{if $p$ is odd} 
\end{cases}
\]
and is generated by a rotation of order $p$ together with a reflection
(with fixed plane orthogonal to the rotation axis).


\section{Symmetry groups and  Lagrangean action}
\label{sec:symm}
Let $\XX$ be \emph{configuration space} of $n$ point particles
in $\RR^3$:
$\XX  = (\RR^3)^n$.
Let $\TT$ be the circle of length $T= |\TT|$.
A function $\TT \to \XX$ is a $T$-periodic path in $\XX$.
By \emph{loops} in $\XX$ we mean the elements 
of  the Sobolev space $\Lambda = H^1(\TT,\XX)$, i.e.
all $L^2$ functions $\TT \to \XX$ with $L^2$-derivative.
The aim is to find 
periodic (in an inertial frame or in a
uniformly rotating frame) orbits 
for the $n$-body problem:
they can be obtained as 
critical points of the Lagrangian action functional
\begin{equation}\label{lag}
\action_\bfomega = 
\int_\TT \left( 
\sum_{i\in \nn}
\frac{m_i}{2} | \dot x_i(t) + \Omega x_i(t)| ^2 
+
\sum_{\stackrel{i<j}{i,j\in \nn}}
m_im_j|x_i(t) - x_j(t)|^ {-\alpha}
\right) \dt,
\end{equation}
where $\Omega$ is the anti-symmetric $3\times 3$ matrix defined by
the relation $\Omega v = \bfomega \times v$ for every
$v\in \RR^3$, with 
the vector $\bfomega\in \RR^3$ representing the rotation
axis of the rotating frame and its norm $|\bfomega|$ 
the angular velocity.
The domain of  the functional $\action_\bfomega$ 
is $\Lambda = H^1(\TT,\XX)$ (of course, allowing a range with infinite
value).
Any collisionless critical point is in fact a $C^2$ solution
of the corresponding 
Euler-Lagrange, or Newton,
equations under a homogeneous graviational potential  
of degree $-\alpha$,
which is periodic in the rotating frame.

Now consider a group $G$ acting orthogonally
on $\TT$, $\RR^3$ and permuting the indices $\nn$.
In other words, consider three homomorphims
$\tau$, $\rho$ and $\sigma$ 
from $G$ to 
$O(\TT)$, $O(3)$ and $\Sigma_n$ respectively.
The group $G$ can be seen as subgroup (possibly $\mod$ 
a normal subgroup) 
of the
direct product $O(\TT) \times O(3) \times \Sigma_n$
under the monomorphism
$\tau\times \rho \times \sigma$,
and the three homomorphisms can be recovered
as projections onto the first, second and third
factor of the direct product.
Given $\rho$ and $\sigma$, 
it is customary to define an action on the configuration
space $\XX$ by  the rule
$(\forall g\in G), x_{\sigma(g)i} = \rho(g) x_i$.
We will denote simply by $gx$ the 
value of $g\cdot x$ under this action of $G$.
In the same way, 
the action of $G$  on $\TT$ and $\XX$ induces 
an action on the functions $\TT \to \XX$ 
by the rule
$(\forall g\in G), x(\tau(g)t) = g x$,
and therefore $\Lambda$ is a $G$-equivariant vector space
(the action of $G$ is orthogonal under the standard
Hilbert metric on $\Lambda$).

\begin{defi}
\label{defi:laggroups}
A subgroup of $O(\TT) \times O(3) \times \Sigma_n$
is termed
\emph{symmetry group}.
It will be termed a 
\emph{symmetry group}
of the Lagrangian  action
functional $\action_\bfomega$ 
if it leaves the value of the action $\action_\bfomega$ 
\ref{lag} 
invariant.
\end{defi}

Note that if $i,j\in \nn$ are indices and 
$gi=j$ for some element $g\in G$,
then it has to be $m_i= m_j$.
More generally, consider the 
decomposition of $\nn$ into (transitive) $G$-orbits,
also known as 
\emph{transitive decomposition}. Indices in the 
same $G$-orbit must share the value of the mass and,
furthermore, the trasitive decomposition yields 
an orthogonal splitting of the configuration space:
\begin{equation}\label{transdec}
\XX  = (\XX_1+\XX_{g1}+\dots) \oplus (\XX_2 + \XX_{g2} + \dots) \oplus \dots,
\end{equation}
where each $\XX_j$ is a copy of $\RR^3$ and each summand  grouped
by brackets is given by 
a transitive $G$-orbit in $\nn$.
This transitive decomposition is nothing but the standard decomposition
of a permutation representation in the Burnside ring $A(G)$.
\begin{defi}
\label{defi:rotaxis}
Consider a symmetry group $G$.
A vector $v\in \RR^3$ is a \emph{rotation axis}
for $G$ if  $(\forall g\in G) gv \in  \{\pm v\}$
(that is, the line $\langle v\rangle \subset \RR^3$ is $G$-invariant)
and the orientation $G$-representation on the time circle
(that is, $\det(\tau)$)
coincides with the orientation 
representation on the orthogonal plane of $v$
(that is, $\det(\rho)\det(v)$).
\end{defi}

We recall from \cite{dd} (proposition 2.15) that
if $\bfomega$ is a rotation axis for a symmetry group
$G$ (and the values of the 
masses are compatible with the transitive
decomposition \ref{transdec}) then $G$
is a symmetry group of the action functional $\action_\bfomega$.
The converse holds, after a straightforward proof, 
for linear or orthogonal actions.
\begin{lemma}
In case the group has a rotation axis it is termed
group \emph{of type R}. 
If the symmetry group $G$ is not of type R, then all $G$-equivariant
loops have zero angular momentum.
\end{lemma}
\begin{proof}
The proof an analogous 
proposition  for $3$ bodies can be found in 
\cite{dd}, proposition 4.2; the details are given for $3$ bodies,
but it can be trivially generalized to the case of $n$ bodies:
if $J$ denotes the angular momentum of the
$G$-equivariant path $x(t)$, 
for every $g\in G$ the formula 
\[
J(gt) = 
\det(\rho(g)) \det( \tau(g)) \rho(g) J(t)
\]
holds,
and hence the angular momentum  $J$ (which is constant) 
belongs to the subspace $V$ in  $\RR^3$ fixed by the $G$-representation
$\det(\tau)\det(\rho)\rho $.
But if $V\neq 0$, then 
 there is a non-trivial vector $v\in \RR$ 
with the property that for every $g\in G$, 
$\det(\tau(g))\det(\rho(g))\rho(g) v = v $.
The orientation representation on the plane orthogonal
if $v$ denotes the representation on $\langle v \rangle$
and $\rho_2$ the representation on its orthogonal compmlement,
it follows that
$\det(\tau(g))\det(\rho_2(g)) \det(v)  \det(v)  = 1 $,
and hence that $\det(\tau) = \det(\rho_2)$: 
the direction spanned by $v$ is a rotation axis, which
contradicts the hypothesis.
\end{proof}

Let $\IsoR$ denote the group of 
(affine) isometries of the time line $\RR$, 
generated by translations and reflections. 
For every $T>0$ there is a surjective projection $\IsoR\to O(\TT)$,
where $\TT= \RR/_{T\ZZ}$.
Let $G$ be a symmetry group and $\widetilde G$ its
\emph{cover} in $\IsoR \times O(3) \times \Sigma_n$,
that is the pre-image of $G$ via the projection
\[
\IsoR \times O(3) \times \Sigma_n \to 
O(\TT) \times O(3) \times \Sigma_n. 
\]
It is easy to see that there is a canonical
isomorphism
\[
H^1(\RR,\XX)^{\widetilde G} \cong H^1(\TT,\XX)^G
\]
and hence that we can consider solutions of 
the $n$-body problem which are $\widetilde G$-equivariant loops
instead of the periodic solutions of the $n$-body problem
which are $G$-equivariant.
Assume now that the symmetry group $G$ has a rotating axis,
and therefore that $\action_\bfomega$
is $G$-invariant. For a fixed angular speed $\theta$, 
the equation $x(t) = e^{i \theta t} q(t)$ 
induces
an isomorphism $\theta_* q \mapsto x$
\[
\theta_* \from H^1(\RR,\XX) \to H^1(\RR ,\XX).
\]
The image $\theta_*\left( H^1(\RR,\XX)^{\widetilde G} \right)$
can be seen as 
\[
\theta_*\left( H^1(\RR,\XX)^{\widetilde G} \right) =
H^1(\RR,\XX)^{\widetilde G'} 
\]
for a new symmetry group $\widetilde G'$ (still of type  R)
with the property that 
if $g$ is a time translation, then $\rho(g)$ is trivial.
Since the following diagram commutes 
\[
\xymatrix{%
 H^1(\RR,\XX)^{\widetilde G} \ar[rr]^\approx 
\ar@/_/[dr]_{\action_\bfomega}  
 &  & H^1(\RR,\XX)^{\widetilde G'} 
\ar@/^/[dl]^{\action_{\bfomega'}} 
 \\
  &  
\protect\RR  & \\
}
\]
(where $\bfomega'$ is 
chosen as suggested above)
a suitable change of angular speed allows 
one to reduce the size of the symmetry group
$G$  and assume that 
(if it is of type R, of course) 
$\ker (\det \tau) \subset \ker \tau \cup \ker \rho $.

\begin{defi}
\label{defi:long}
Consider the following terms.
A symmetry group $G$ is:
\begin{description}
\item[bound to collision:] if every $G$-equivariant loop
has collisions;
\item[homographic:] if every $G$-equivariant loop is homographic, 
i.e. constant up to Euclidean similarities;
\item[transitive:] if the prermutation action
of $G$ on the index set is transitive;
\item[fully uncoercive:] if for every possible rotation vector
$\bfomega$ the corresponding action functional 
$\action_\bfomega$ is coercive if restricted to the 
space of $G$-equivariant loops $\Lambda^G$.
\end{description}
\end{defi}
\begin{defi}
\label{defi:core}
The kernel $\ker \tau$ is termed the \emph{core} of the symmetry group.
\end{defi}
\section{Transitive groups and transitive decomposition}
\label{sec:trans}
Consider a symmetry 
group $G\subset O(\TT) \times O(3) \times \Sigma_n$
and its \emph{cover}
$\widetilde G \subset \IsoR \times O(3) \times \Sigma_n$
(which anyhow 
is a discrete group acting on the time line $\RR$ 
as time--shifts and time--reflections); the kernel of the projection
$p\from \widetilde G \to G$ is a free abelian group of rank $1$.
By composition with the projection $p\from \tilde G \to G$ 
it is possible to define
the homomorphisms 
$\tilde \tau = \tau p$,
$\tilde \rho = \rho p$ and
$\tilde \sigma = \sigma p$.
Let us note that the diagram
\[
\xymatrix{%
{\widetilde{G}} \ar@{ >->}[r]^-{\tilde \tau \times \tilde \rho} 
\ar@{->>}[d] & {\IsoR \times O(3) }
\ar@{->>}[d] 
\\
G \ar@{ >->}[r]_-{\tau \times \rho} & {O(\TT) \times O(3)} \\
}
\]
commutes: the horizontal arrows  are monomorphisms and 
the vertical arrows are epimorphisms.
The projection onto the second factor  of 
$\tG \subset \IsoR \times O(3)$ has as image a finite space point group 
$F\subset O(3)$.
Now consider the core of $G$, 
$\ker \tau  = K\subset G$. Its pre-image $\widetilde K = p^{-1}(K) =  \ker\tilde\tau \subset \tG$
is isomorphic to $K$ via $p$  
and 
the restriction of $\rho$ to $K$ is a monomorphism.
With an abuse of notation
we can identify 
$K \cong \widetilde K  \subset \widetilde G \subset G\subset \IsoR \times O(3)$  
with its image in $O(3)$ under $\tilde \rho$
Consider the normalizer $N_{O(3)}K$ of $K$ in
$O(3)$. Since the core $K$ is normal in $G$, 
the space group $F$ is a subgroup of the normalizer $N_{O(3)}K$.
If follows that:
\begin{nr}
There is a homomorphism
$\widetilde G / \widetilde K \to G/K \to  W_{O(3)}K$
of $G/K$ with image in the Weyl group of $K$
in $O(3)$.
\end{nr}
Now note that 
$p^{-1}(\ker\det\tau) = \ker \det \tilde \tau$ and 
consider the fact that the quotient
$\det\ker\tilde\tau/K$ 
is isomorphic to  $\ker p \cong \ZZ$.
In fact, $G/K$ is projected onto a subgroup of 
$O(\TT) \times W_{O(3)}K$,
while 
$\widetilde G/ \widetilde K$ is projected
onto a subgroup of $\IsoR \times W_{O(3)}K$.

\begin{defi}\label{defi:Krh}
Let $G$ be a symmetry group. Then define
\begin{enumerate}
\item $K=\ker\tau$, 
\item 
$[r]\in W_{O(3)}K$ as the image in 
the Weyl group of 
the generator mod $K$
of $\ker\det\tau \subset G/K$ 
(corresponding to the time--shift  with minimal angle).
If $\ker\det\tau = K $, then $[r]=1$.   
\item
$[h]\in W_{O(3)}K$ as the image in the Weyl group
of one of the time-reflections mod $K$ in $G/K$,
in the cases such an element exists. Otherwise it is not defined.
\end{enumerate}
In short, the triple 
\[
(K,[r],[h])
\]
is said the $\Krh$ data of $G$.
\footnote{%
For the sake of simplicity we omit such square brackets when unnecessary.}
\end{defi}
Now assume that the action of $G$  on the index set
is transitive. Then the following easy results hold.
\begin{nr}\label{nr:isotropy}
For all $i\in \nn$, the isotropy subgroups
$H_i = \{g\in G \st gi=i \}$ are mutually conjugated in $G$.
Left multiplication by elements in $G$ yields 
a bijection $G/H_1 \cong \nn$ 
of the index set $\nn$ and the set of left cosets $G/H_1$,
which is $G$-equivariant (that is, a $G$-bijection).
\end{nr}

We can here define the last piece of information needed for 
a classification of the symmetry groups: let $H_1$ denote
one of  the isotropy subgroups defined above in \ref{nr:isotropy}.

\begin{nr}
\label{nr:threecases}
Assume that $\ker\tau \neq 1$. Then 
one (and only one) of the following cases can occur:
\[
\begin{cases}
\ker\tau \cap H_1 = 1 \\
\ker\tau \cap H_1 = \langle \text{reflection along a plane} \rangle \\
\ker\tau \cap H_1 = \ker\tau. \\
\end{cases}
\]
\end{nr}
\begin{proof}
We tacitely assumed that the core 
$\ker\tau$ is not a reflection along a plane,
since otherwise the problem would be a planar $n$-body problem
or bound to collisions.
Furthermore, 
since we assume $\ker\tau \neq 1$, the \emph{only one} part is trivial.
Suppose on the other hand 
that $\ker\tau \cap H_1 \neq 1$. Let $E\subsetneq \RR^3$ be the 
linear fixed  subspace
\[
 E = \left( \RR^3\right)^{\ker\tau \cap H_1}.
\]
Let $K$ be as above $K=\ker\tau$.
The configuration space $\XX$ can be seen as the space
of maps $G/H_1\to\RR^3$, where $G/H_1$ is seen as a $G$-set
with $[G:H_1]$ elements and $\RR^3$ is of course a $G$-space via $\rho$.
The action on $\XX$ (as space of maps) is the diagonal action,
and configurations in $\XX^K$ correspond to 
$K$-equivariant maps $G/H_1\to \RR^3$.
Now, 
the number of $K$-orbits in $G/H_1$ is also
the number of the double cosets $K\backslash G/H_1$;
since $K$ is normal in $G$, it coincides with the 
number of $H_1$-orbits in $G/K$, which is 
$[G:H_1K]$.
Any $K$-map $x\from G/H_1 \to \RR^3$  
(i.e. an element in $\XX^K$)
can therefore be decomposed into a sum of $[G:H_1K]$ disjoint parts 
(more precisely,
its domain can be) which are the $K$-orbits in $G/H_1$. 
Each map defined on a $K$-orbit is conjugated
via an element of $G$  
to a $K$-map of type
\[
K/(K\cap H_1) \to (\RR^3)^{|K\cap H_1|} = E,
\]
(thus yielding $[K:K\cap H_1]$ particles
in $E$).
The space $\XX^K$ is isomorphic 
to a direct sum of $[G:H_1K]$ copies of $E$,
over which the action of $G$ acts via conjugation
(actually, it is the  induced/inflated module).

Now, consider the for hypothesis that  $K\neq K\cap H_1$.
The dimension $\dim E$  can be $0$, $1$ or $2$ (it cannot be $3$ since
by assumption $K\cap H_1 \neq 1$ and the action of $K$ on $\RR^3$ is
faithful).
If it is $0$, then $K=K\cap H_1$ (since otherwise at each time a collision
would occur.
If $\dim E = 1$ (and $K\neq K\cap H_1$, assumption above)
then it is easy to show that either the group $G$ is fully uncoercive
or it is bound to collisions:
in fact for one-dimensional $E$ there cannot exist rotation axes,
and any symmetry element yielding coercivity would make
the group bound to collisions (the complementary of the collision
set in $E$ is not connected).
It is left the case $\dim E=2$, i.e. where $K\cap H_1$ is 
the group generated by a single plane reflection.
If the plane $r$  fixed by $K\cap H_1$ is $K$-invariant (that is,
$Kr=r$), then $K\cap H_1$ is normal in $K$ and the $K$-representation
given by $\rho$ is one of the following:
\begin{enumerate}[\itshape 1)] 
\item $C_{ph}$ with $p\geq 1$ (the group generated by the reflection
around the plane $r$ and $p$ rotations orthogonal to $r$), 
\item $I\times D_p$ with $p\geq 2$ even (the Coxeter group generated by
the reflection around $r$ and $p$ ``vertical'' plane reflections),
\item $D_pC_p$ with $p\geq 2$ (the Coxeter group generated by $p$ plane reflections) 
and 
\item $D_{2p}D_p$ with $p\geq 1$ (generated by $D_p$ and $-\zeta_{2p}$: 
 it is a Coxeter group for $p$ odd).
\end{enumerate}
Cases  2), 3)  and 4)  do not possibly have rotation axes,
and a symmetry group extending $K$ and not coercive would be
fully uncoercive.  
Since the bodies are constrained to belong to $r$ ($H_1$ is the isotropy
of the permutation action) and the singular set of
$r$ cuts $r$ into different components, a symmetry group extending
such $K$ cannot be coercive without being
bound to collisions.
Case 1) is of a different type: the $p=[K:K\cap H_1]$ bodies are
constrained to be vertices of a regular $k$-agon centered at the origin 
and contained in $r$. The direction orthogonal to $r$ is a rotation axis.
This is 
the case in which  the reflection along a plane yields possible 
periodic orbits.
\end{proof}
Note that 
if $\ker\tau \cap H_1 = 1$, then 
the isotropy $H_1$ is isomorphic to its 
image under $\rho$ 
(after the composition with the projection onto the $K$-quotient)
in the Weyl group  $W_{O(3)}K$.

\begin{defi}
\label{defi:hKrh}
Let $G$ be a symmetry group acting transitively
on the index set and $H_1\subset G$ the isotropy subgroup
with respect to the index permutation action,
defined in \ref{nr:isotropy}. If $\widetilde{H_1}$ is 
the cover of $H$ 
(that is, $p^{-1}H_1 \subset \tG \subset \Iso(\RR) \times O(3)$)
Consider the following pieces of data:
\begin{enumerate}
\item 
$\widehat K
= K \cap \widetilde H_1 
\subset \tG$
(note that $\widehat K \cong \ker\tau\cap  H_1$);
that is, elements of $\widehat K$ are the elements in $\widetilde{H_1}$
fixing the time.
\item 
Consider the image of $H$ in $G/K$. Its intersection 
with the cyclic group $\ker\tau/K$ is a cyclic group with 
a distinguished non-trivial generator (if not trivial),
say $r \mod K$. One of its pre-images
 $p^{-1}r$ in 
$\IsoR \times N_{O(3)}K \subset \IsoR \times O(3)$
is an element 
$(k,\widehat r)$
in $\widetilde H_1 \cap \ker\det\widetilde \tau$,
defined up to multiplication with elements in $\widehat K$.
Without loss of generality one can assume $k$ to be an integer.
\item
If the set $\widetilde H \cap (\widetilde G \minus \ker\det\widetilde \tau)$
is non-empty, then 
let $\widehat h$ be the projection  in 
$N_{O(3)}K$ of one of its elements.
\end{enumerate}
Then the triple
\[
(\widehat K, (k,\widehat r), \widehat h)
\]
is said to be the $\hKrh$ data of $G$. 
\end{defi}
\begin{propo}
\label{propo:struct}
Let $G$ be a transitive symmetry group and 
$\begin{pmatrix}
\widehat K &   (k, \widehat r) &  \widehat h \\
K & [r] & [h] \\
\end{pmatrix}$ 
the matrix with as first row that $\hKrh$ 
data defined above in \ref{defi:hKrh} 
and as second row the $\Krh$ data defined in \ref{defi:Krh}.
Then the cover $\tG$ of $G$ is, up to conjugacy, defined in 
$\IsoR \times O(3)$  by the $\Krh$ data.
The cover of the isotropy  
$\widetilde H_1$ is defined by the $\hKrh$ data in the first row.
Its permutation representation on indices can be deduced  
by considering that $G/H_1 \cong \widetilde G / \widetilde H_1$. 
\end{propo}

This allows to properly define a decomposition of the 
permutation action of $G$ into transitive components,
using the following sum.

\begin{defi}
Let $G_1$ and $G_2$ two groups with the same $\Krh$ data.
Then the disjoint sum $G_1 + G_2$ is defined 
as follows:
the covers $\widetilde{G_1}$ and 
$\widetilde{G_2}$ are isomorphic and generated in $\IsoR \times O(3)$ 
by the 
(common) $\Krh$ data.
The action of such a resulting $\tG$ 
on the index set  can be defined by taking the disjoint
union of the $\tG$-sets $\tG/\widetilde H_1 + \tG/\widetilde H_2$.
Now to find the projection $p\from \tG \to G$,
it suffices to consider the least common multiplier 
of the integers $k_1$, $k_2$, $|r_1 \mod \widehat K_1|$
and $|r_2 \mod \widehat K_2|$.
\end{defi}

\section{Local variations and averaging techniques 
over equivariant spheres}
\label{sec:averaging}
The sum defined in the previous section allows one to build 
and generate all symmetry groups using their transitive components.
Now we come to the problem of collisions.
Let $\tG$ be (the cover of) a symmetry group and 
$x = x(t) \in \Lambda = H^1(\RR,\XX)^\tG$
a local minimizer. Assume that  at time $t=0\in \RR$
the trajectory $x(t)$ collides, and all bodies 
in a cluster $\kk \subset \nn$ collide (which means that other bodies
might collide, but not with bodies in $\kk$). 
In sections 7--9 of \cite{FT2003} the blow-up and the 
averaging technique are developed for equivariant trajectories; we 
refer to it for details.
We now extend the range of applicability of the averaging techique to
symmetry groups that do not need to fulfill the rotating circle property.
The blow-up of $x(t)$ centered at $0$
$\bar q$ is a local minimizer of the Lagrangian action $\action$
restricted to the space
$H^1(\RR,\XX_\kk)^{G_*}$,
where $G_*$ is the restriction of $\tG$ to the subgroup 
\[
\IsostarR \times O(3) \times \Sigma_\kk 
\subset
\IsoR \times O(3) \times \Sigma_n,
\]
$\IsostarR$ is the group of order $2$ 
consisting in the isometries of $\RR$ fixing
$0$ and $\Sigma_\kk$ the permutation group on 
the indices in $\kk \subset \nn$.
In other words, $G_*$ is the 
symmetry group with $\Krh$ data
$(K,1,h)$, the  transitive 
decomposition obtained restricting the permutation
action to the colliding particles in $\kk$
and $\XX_\kk$ denotes the configuration space
of the particles in $\kk$.
Another way to define $G_*$ is to consider it 
as the subgroup of all elements in $\tG$ (or, equivalently, $G$) 
fixing the colliding time (that is, its isotropy
subgroup, a maximal isotropy subgroup).

Now it comes to the standard variation.
Let  us define
\[
S(s, \delta) = 
\int_0^\infty
\left[
\dfrac{1}{|t^{2/(2+\alpha)}s + \delta|^\alpha} 
-
 \dfrac{1}{|t^{2/(2+\alpha)}s|^\alpha}
\right]   \dt 
\] 
The following lemma  follows from  section 9 of \cite{FT2003}.
\begin{nr} 
\label{nr:standardvariation}
Let $\bar q(t)$ be a colliding blow-up trajectory 
and $\bar s$ the limiting central configuration in $\XX_\kk$.
If there exists a symmetric  
configuration $\delta \in \XX_\kk^{G_*}$ (that is, $\delta$ is 
fixed by the 
isotropy $G_*$) such that for every 
$i,j\in \kk$
\[
S(\bar s_i - \bar s_j, \delta_i - \delta_j) \leq 0
\]
and for at least a pair of indices the inequality is strict, 
then 
the colliding blow-up trajectory  $\bar q(t)$ 
is not a minimizer.
\end{nr}

Now we consider three different procedures that can be 
used to find such a $\delta$. A symmetric variation $\delta$
that let the action functional $\action$ decrease on the 
standard variation is called \emph{V-variation}.

The following proposition is 
contained in the proof of theorem (10.10) of  \cite{FT2003}.
\begin{nr}
\label{nr:proc1}
If $G_*$ acts trivially on $\kk$, then
a V-variation always exists. 
\end{nr}

We recall that we say that a circle $\SS\subset \RR^3$ (with center
in the origin $0$) is 
called \emph{rotating} under a group $G_*$ for an index when
it is $G_*$-invariant and 
$\SS \subset \left(\RR^3\right)^{H_i}$,
where $H_i\subset G_*$ is the isotropy of $i$ with respect to the 
the permutation action of $G_*$ on the index set, via $\sigma$.
Now, proposition (9.8) of \cite{FT2003} can be re-phrased as follows.
\begin{nr}
\label{nr:proc2}
If there is an index $i\in \kk$ and a circle
$\SS \subset \RR^3$ which is rotating under $G_*$ for the 
index $i$, then 
the average
\[
\int_{\delta \in \iota_i\SS}
\int_{\delta \in \iota_i\SS}
\sum_{j\neq i}
S( \bar s_i - \bar s_j, \delta_i - \delta_j)
< 0,
 \]
is strictly negative,
where $\iota_i \SS \subset \XX_\kk$ is the image of the rotating circle $\SS$ 
under the inclusion $\iota_i$ 
defined 
as the inclusions in the proof of \ref{nr:threecases}.
In other words, if there is a rotating circle under $H_i$ then
by averaging it is possible to find a V-variation.
\end{nr}

Note that \ref{nr:proc2} holds true if and only if 
the hypothesis of the claim 
is true for even just one of the transitive components 
of the index set.
In other words, a V-variation obtained
by averaging over a circle
exists if and only if it is possible to obtain a V-variation
by averagin over a circle only
in one of the transitive components in which $G_*$ 
can be subdivided.
The same will be true for the next proposition,
which is a new generalization of the rotating circle property.
\begin{nr}
\label{nr:proc3}
Let $G_*$ be the symmetry group of a blow up solution $\bar q$ as above. 
If $\det\rho(G_*)=1$ (i.e. $G_*$ acts orientation-preserving on 
the space $\RR^3$) and 
for one of the indices $i\in \kk$ the permutation isotropy $H_i$
(restricted to $G_*$) is trivial,  
then there exists a V-variation, obtained 
by averaging over a $2$-sphere. 
\end{nr}
\begin{proof}
Let $S^2\subset \RR^3$ be a $2$-sphere
centered in $0$. 
If $H_i = H_i \cap G_* = 1$,
then the space $E=\left( \RR^3\right)^{G_* \cap H}$ 
defined as in the proof of \ref{nr:threecases} 
is equal to $\RR^3$ and it contains the sphere $S^2$.
As exaplained in the same proof, 
the fixed configuration space $\XX^{G_*}$
can be decomposed into a sum of some copies of 
$E$ (exactly $|G_*|$, since the isotropy is trivial)
and a remainder (which depends on the indices which are
not in the same homogeneous part of $i$):
there is hence an embedding
$\iota_i\from S^2 \to \XX_\kk$ defined by the group action.
Now, all elements of $G_*$ by hypotheses act by rotations
on $\RR^3$.
Now consider the average
\[
A = 
\int_{\delta \in \iota_i S^2}
\sum_{i<j}
S ( \bar s_i - \bar s_j , \delta_i - \delta_j ) 
\]
Each term in the sum is equal to 
the sum of terms like
\[
A_g = \int_{\delta_i \in  S^2}  
S ( \bar s_i - \bar s_j , (1-g)\delta_i ) 
\]
where $g$ range in $G_*$. But since $g$ acts as rotation in $\RR^3$,
$(1-g)$ is a projection onto a plane composed with a dilation:
for each $g\in G_*$ therefore  
there is a positive constant $c_g>0$ such that 
 \[
A_g = c_g 
\int_{\delta_i \in  \SS} 
S ( \bar s_i - \bar s_j,  \delta_i )
\]
obtained exactly as in the case of the integration of a disc.
Since such terms  are strictly negative,
the conclusion follows.
\end{proof}

\section{Transitive components of groups with non-colliding minimizers}
\label{sec:proof}

\begin{defi}
\label{todo:assumption}
We say that a group $G$  has property
\ref{todo:assumption} if it is:
\begin{inparaenum}
\item  not bound to collision,
\item  not fully uncoercive, 
\item  not homographic and at last that 
\item\label{prop:last}  for all maximal time-isotropy subgroups $G_* \subset G$ 
at least one of the propositions
\ref{nr:proc1},
\ref{nr:proc2}
or \ref{nr:proc3}
can be applied (that is, either $G_*$ acts trivially on indexes,
or there is a transitive component with a rotating circle or 
$G_*$ acts by rotations).
\end{inparaenum}
\end{defi}

If the group is not fully uncoercive, then possibly considering
a non-zero angular velocity $\omega$ it is possible to show 
that local minima always exist. We exclude the groups bound to collisions
and homographic simply because we are looking for collisionless 
solutions which are not homographic.
Now, if furthermore property \ref{prop:last}
(which can be easily tested only on the transitive components, 
as noted above) holds, the existence 
of a V-variation implies that all local minimizers are collisionless,
which is our goal.
We start by considering the possible cores for $G$ (not considering 
at the moment the permutations on the indices).

All finite subgroups of $SO(3)$ listed in  table \ref{tab:finitegroups1} 
(and the trivial
group, not listed)  
can be cores by \ref{nr:proc3}, as far as the isotropy ($\widehat K$
in the $\hKrh$ data of the corresponding component)  of one of
the indices is trivial.
Then, of the groups of table \ref{tab:finitegroups2},
the central inversion group $I$ and the central prism/antiprysm 
$I\times C_p$ group have a rotating circle  
and can be considered.
The groups $I\times D_p$, with $p\geq 2$ are generated by plane
reflections for $p$ even and do not contain rotating circles:
the action restricted to  invariant planes is never consisting of 
rotations. The only possible hypothesis for the existence of a V-variation
is the triviality of the permutation action: but the subspace of $\RR^3$ 
fixed is $0$, and hence with more than one particle the group
would be bound to collisions.
The remaining groups $I\times T$,
$I\times O$ and $I\times Y$ of the table act on $\RR^3$ without invariant planes
(the representation is irreducible)
and hence they must be excluded. 
The same is true for the full tetrahedron group  $OT$  of table
\ref{tab:finitegroups3}.
Of the three remaining groups in this table,
the prism/antisprims group $C_{2p}C_p$ clearly has a rotating circle
and must be added to the list.
The  groups $D_pC_p$ (the $p$-gonal planes reflection 
group) and $D_{2p}D_p$ (for $p\geq 2$) do not have rotating circles
and have reflections: not only none of the \ref{nr:proc1},
\ref{nr:proc2} and \ref{nr:proc3} can be applied: all the symmetry
groups having this core would result to be bound to collisions
or fully uncoercive.

\begin{nr}\label{nr:fulllist}
The groups satisfying 
\ref{todo:assumption} are  the following:
\begin{inparaenum}[\itshape 1)]
\item $C_p$ {(for $p\geq 1$)},
\item $I\times C_p$ {(for $p\geq 1$)},
\item $C_{2p}C_p$ {(for $p\geq 2$)},
\item $D_p$ {(for $p\geq 2$)},
\item $T$,
\item $O$,
\item $Y$.
\end{inparaenum}
\end{nr}

Of course for the same reason this  is also the list of 
projections on $O(3)$ of the (possible) maximal time-isotropy groups
and of the cores.
The same argument yields the proof of 
Theorem \ref{theo:main}.
Now we consider the extensions (of index $2$) of such cores
as possible time-istropy for times fixed by reflections. 
The method used for obtaining the existence of V-variations
sets constraints on the type of admissible extensions: 
a group with V-variations only by avaraging on spheres and 
without rotating circles cannot be extended other than in 
$SO(3)$, as in the case of the last four items in the list.
On the other hand, the two prism/antiprysm family of groups
$I\times C_p$ and $C_{2p}C_p$ have a rotating plane
but they are not orientation-preserving: hence
can be extended without restrictions on the orientation
but keeping the rotating plane.
Also, we need to rule out from the list of possible cores 
the groups that  do not occur as cores 
of symmetry groups not bound to collisions and not fully uncoercive.
By the same argument used in the proof of \ref{nr:threecases},
if $\RR^3$ is disconnected by the collision subspaces, then 
it is not possible to assume coercivity and being collisionless. 
Hence those groups with fixed planes must be eliminated:
of the two families $I\times C_p$ and 
$C_{2p}C_p$ only the anti-prism family of groups
$P'_{2p}$ survives, with normalizer $I\times C_{2p} (= C_{2ph})$.
Simple geometric and algebraic considerations lead us to the 
following conclusions:
\begin{nr}\label{nr:fulllist2}
The 
index $2$ extensions 
satisfying 
\ref{todo:assumption} 
of  
cores  
satisfying \ref{todo:assumption}   
are  the following:
\begin{inparaenum}[\itshape 1)]
\item $C_1$: $I$, $C_2C_1$, $C_2$.
\item $C_p$ {(for $p\geq 2$)}: $C_{2p}$, $D_p$, $I\times C_p$, 
$C_{2p}C_p$. 
\item $P'_{2p}$ {(for $p\geq 1$)}: $I\times C_{2p}$,
\item $D_p$ {(for $p\geq 2$)}: $D_{2p}$.
\item $T$: $O$.
\item $O$: nothing.
\item $Y$: nothing.
\end{inparaenum}
\end{nr}

Recall that the matrices of $\Krh$ and $\hKrh$ data  are
$\begin{pmatrix}
\widehat K &   (k, \widehat r)  \\
K & [r]  \\
\end{pmatrix}$  (for the cyclic type)
or 
$\begin{pmatrix}
\widehat K &   (k, \widehat r) &  \widehat h \\
K & [r] & [h] \\
\end{pmatrix}$ 
(for brake or dihedral type),
as defined in  \ref{propo:struct}. 
\subsection*{Trivial core}
Let us now consider the simpler case of trivial core. By definition 
$K=1$ and $\widehat K=1$, which implies 
$\ze \cong \tG =\langle (1,r)\rangle 
\subset \Iso(\RR) \times O(3)$. About the pair $(k,\widehat r)$
generating the cover $\widetilde H_1$ of the permutation  isotropy, 
it must be a power of the generator $(1,r)$ and hence
of the form $(k,r^k)$.

If the action is of cyclic type, then the  $\Krh$ can be written as
$\begin{pmatrix}
1 &   (k, r^k)  \\
1 & r  \\
\end{pmatrix}$, where 
up to rotating frames $r$ can be chosen with order at most $2$ 
(it is not difficult to see that 
every cyclic symmetry group is of type R).
Since if $r=1$, then it must be $\widehat =1$, 
we have for every $k\geq 1$
the \emph{choreographic} symmetry
\[
\begin{pmatrix}
1 &   (k, 1)  \\
1 & 1  \\
\end{pmatrix},
\]
which acts transitively on the set of $k$ bodies.
Of course, the constraints can be written also as the better known form
$x_1(t+i) =  x_i(t)$ for $i=1,\dots, k$ 
for $k$-periodic loops.

If $r$ is the reflection $-\zeta_2$, then 
the $\Krh$ is 
\[
\begin{pmatrix}
1 &   (k, (-\zeta_2)^k )  \\
1 &  -\zeta_2  \\
\end{pmatrix},
\]
which acts again on set of  $k$ indices, but with a resulting cyclic group
$G$ with $2k$ elements.
Any other choice of $r$ would give rise to one of these groups,
up to a change of rotating frame.

Following the same argument as in section 6 of \cite{dd},
one can see that 
the $\Krh$  for a dihedral group of type R can be chosen of the
following forms (for $h_1$ and $h_2$ integers):
\[
\begin{pmatrix}
1 &   (k, 1) &  *  \\
1 & 1 & (-1)^h_1 \zeta_2^{h_2}    \\
\end{pmatrix}
\text{\ \ or \ \ } 
\begin{pmatrix}
1 &   (k, (-\zeta_2)^k )  & *  \\
1 &  -\zeta_2  &  (-1)^h_1 {\zeta_2}^{h_2} \\
\end{pmatrix}.
\]
Groups not of type R  can be found in an similar fashion.

\section{A few examples}
\label{sec:examples}
\begin{example}
\label{ex:icosahedral}
Consider the icosahedral group $Y$ of order $60$.
The group $G$ with $\Krh$ data 
\[
\begin{pmatrix}
1 &   (1,-1)   \\
Y &  -1  \\
\end{pmatrix}
\]
is isomorphic to the direct product $I \times Y$ of order $120$,
and acts on  the euclidean space $\RR^3$ as the full icosahedron group.
The action on $\TT$ is cyclic and given by the fact
that $\ker \tau = 1 \times Y$. 
The isotropy is generated by the central inversion $-1$, 
and hence the set of bodies is $G/I \cong Y$.
Thus at any time $t$ the $60$ point particles are  constrained 
to be a $Y$-orbit in $\RR^3$ (which does not mean they are vertices
of a icosahedron, simply that the configuration is $Y$-equivariant).
After half period every body is in the antipodal position:
$x_i ( t + T/2) = - x_i$ (in other words, the group
contains the \emph{anti-symmetry},
also known as \emph{Italian symmetry} -- see \cite{ambrosetticoti,ambrocoti,chenICM,chenkyoto}.
Of course, the group $Y$ is just an example: one can choose also
the tetrahedral group $T$ or the octahedral $O$ and obtain
anti-symmetric orbits for 12 (tetrahedral) or 24 (octahedral) bodies,
as depicted in figure \ref{fig:icosahedral}. The action
is by its definition transitive and coercive;
local minimizers are collisionless since the maximal $\TT$-isotropy
group acts as a subgroup of $SO(3)$ (i.e. orientation--preserving).
\begin{figure}
\centering
\hfill
\includegraphics[width=0.3\textwidth]{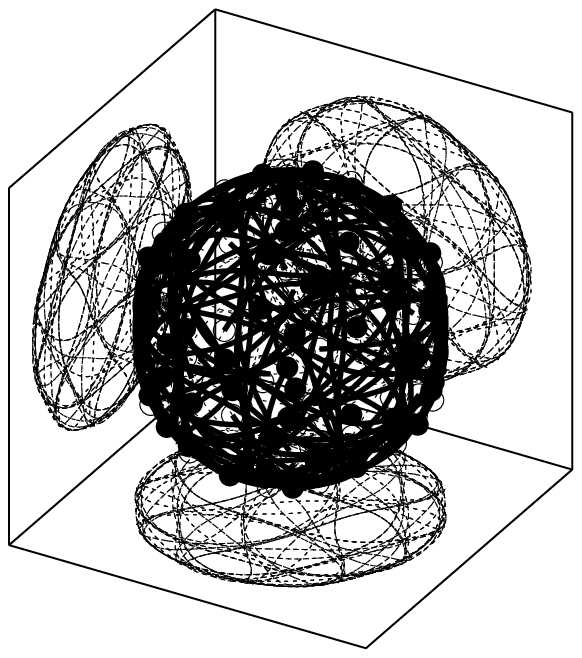}
\hfill
\includegraphics[width=0.3\textwidth]{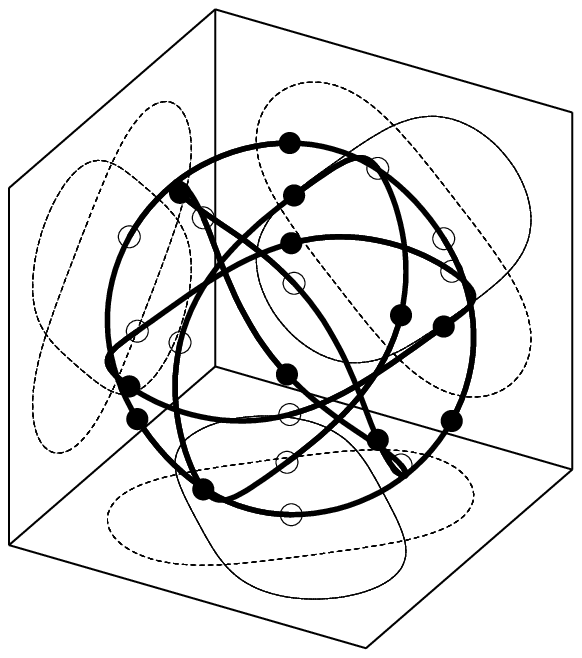}
\hfill
\includegraphics[width=0.3\textwidth]{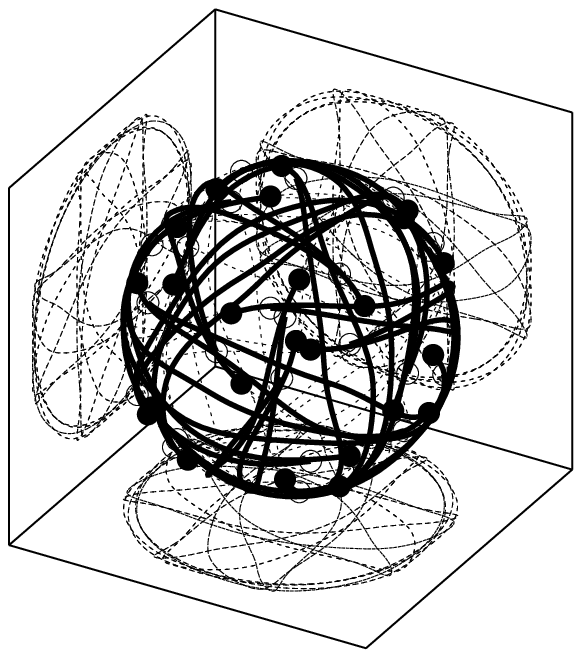}
\hspace{\fill}
\caption{60-icosahedral $Y$, 12-tetrahedral $T$ and 24-octahedral $O$ 
periodic minimizers (chiral)}
\label{fig:icosahedral}
\end{figure}
\end{example}
\begin{example}
\label{ex:klein}
Let $G$  be the group with $\Krh$ data 
\[
\begin{pmatrix}
1 &   (1,-1)   \\
D_k &  -1  \\
\end{pmatrix},
\]
where $D_k$ is the rotation dihedral group of order $2k$.
As in the previous example, the action is such that the action 
functional is 
coercive and its local minima collisionless. At every time instant
the bodies are $D_k$-equivariant in $\RR^k$
and the anti-symmetry holds. Approximations of minima can be seen
in figure \ref{fig:klein}.
\begin{figure}
\centering
\hfill
\includegraphics[width=0.3\textwidth]{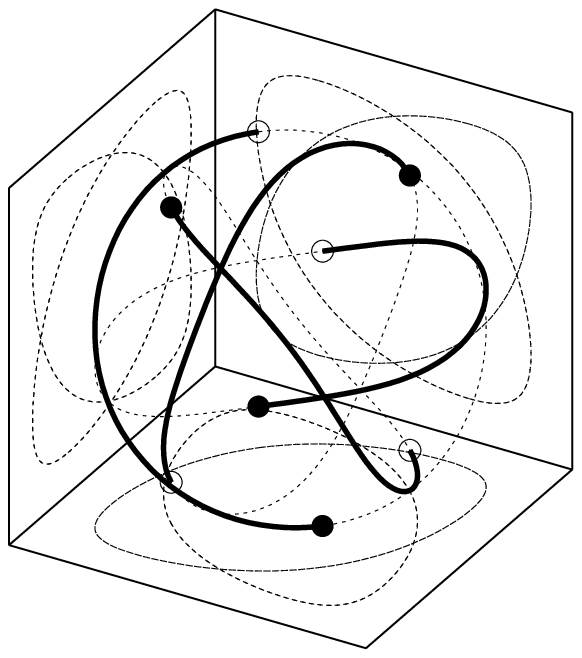}
\hfill
\includegraphics[width=0.3\textwidth]{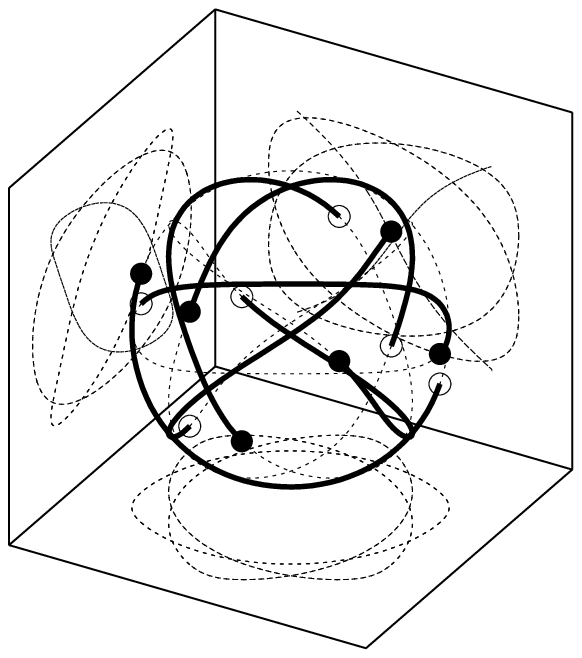}
\hspace{\fill}
\caption{4-dihedral $D_2$  and 6-dihedral $D_6$ symmetric periodic minimizers}
\label{fig:klein}
\end{figure}
\end{example}
\begin{example}
\label{ex:threelag}
To illustrate the case of non-transitive symmetry group,
consider the following (cyclic) $\Krh$ data:
\[
\begin{pmatrix}
1 &   (3,-1)   \\
1 &  -1  \\
\end{pmatrix},
\]
which yields a group of order $6$ acting cyclically 
on $3$ bodies, and with the
antipodal map on $\RR^3$. Since $\ker\tau$ is trivial and   
the group is of cyclic type, local minima are collisionless.
Now, by adding $k$ copies of such group one obtains 
a symmetry group having $k$ copies of it as its transitive components,
where still local minimizers are collisionless and the 
restricted functional is coercive.
Some possible minima can be found in figure \ref{fig:threelag},
for $k=3,4$.
\begin{figure}
\centering
\hfill
\includegraphics[width=0.3\textwidth]{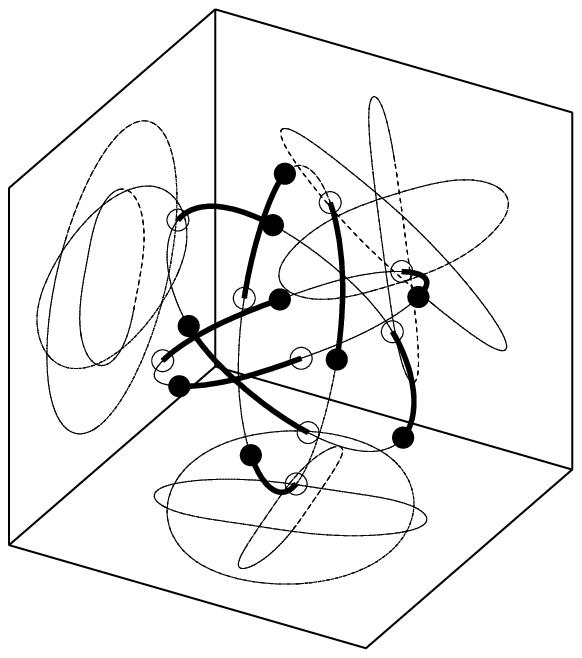}
\hfill
\includegraphics[width=0.3\textwidth]{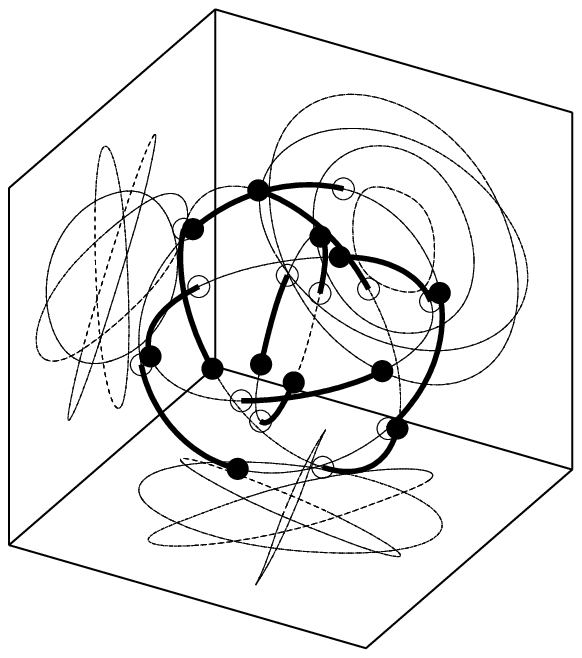}
\hspace{\fill}
\caption{9 and 12 bodies in anti-choreographic constraints grouped by 3}
\label{fig:threelag}
\end{figure}
\end{example}

\begin{remark}
The planar case can be dealt exactly as we did for the spatial case, with
a significative simplification: only when the permutation action is trivial
or there exists a rotating circle (that is, under these hypotheses 
the maximal $\TT$-isotropy
group of all possible colliding times has transitive components which 
act on the position space as rotations).
A transitive decomposition of such planar symmetry group, also,
is much simpler since the core has to be a (regular polygon) cyclic group.
Nevertheless, also in the planar case many examples can be built
using these simple building block. It is still an open problem
whether there are symmetry groups not bound to collisions 
with (local or global?)  minimizers which are colliding trajectories.
It has been proved in \cite{zz} that it cannot happen for $n=3$,
but to the author's  knowledge there is not yet  a general result.
\end{remark}

\nocite{chenICM,MR2146342,MR2032484,MR2020261,MR1688886,monchen,ambrosetticoti,zz}
\def\cfudot#1{\ifmmode\setbox7\hbox{$\accent"5E#1$}\else
  \setbox7\hbox{\accent"5E#1}\penalty 10000\relax\fi\raise 1\ht7
  \hbox{\raise.1ex\hbox to 1\wd7{\hss.\hss}}\penalty 10000 \hskip-1\wd7\penalty
  10000\box7} \def\cprime{$'$} \def\cprime{$'$} \def\cprime{$'$}
  \def\cprime{$'$} \def\cprime{$'$} \def\cprime{$'$} \def\cprime{$'$}
  \def\cprime{$'$} \def\cprime{$'$} \def\cprime{$'$} \def\cprime{$'$}
  \def\cprime{$'$} \def\cprime{$'$} \def\cprime{$'$} \def\cprime{$'$}
  \def\cprime{$'$} \def\polhk#1{\setbox0=\hbox{#1}{\ooalign{\hidewidth
  \lower1.5ex\hbox{`}\hidewidth\crcr\unhbox0}}}
  \def\polhk#1{\setbox0=\hbox{#1}{\ooalign{\hidewidth
  \lower1.5ex\hbox{`}\hidewidth\crcr\unhbox0}}} \def\cprime{$'$}
  \def\cprime{$'$} \def\cprime{$'$}


\end{document}